\documentclass[11pt]{amsart}
\input rlepsf

\newcommand{\bfz}{{\mathbb {Z}}}
\newcommand{\bfq}{{\mathbb {Q}}}
\newcommand{\bfr}{{\mathbb {R}}}
\newcommand{\Li}{\mathbb {L}}
\newcommand{\lk}{{\ell k}}
\newcommand{\rot}{\mbox{\rm rot}}

\begin{document}

\title{Lutz twist and contact surgery}
\author{Fan Ding}
\author{Hansj\"org Geiges}
\author{Andr\'{a}s I. Stipsicz}
\address{Department of Mathematics, Peking University, Beijing 100871, 
P.R. China}
\email{dingfan@math.pku.edu.cn}
\address{Mathematisches Institut, Universit\"at zu K\"oln, 
Weyertal 86-90, 50931 K\"oln, Germany}
\email{geiges@math.uni-koeln.de}
\address{A. R\'enyi Institute of Mathematics, Hungarian Academy of Sciences,
Budapest, Re\'altanoda utca 13--15, H-1053, Hungary}
\email{stipsicz@math-inst.hu}
\date{}
\begin{abstract}
For any knot $T$ transverse to a given contact structure on a $3$-manifold,
we exhibit a Legendrian two-component link $\Li =L_1\sqcup L_2$
such that $T$ equals the transverse push-off of
$L_1$ and contact $(+1)$-surgery on $\Li$
has the same effect as a Lutz twist along~$T$.
\end{abstract}
\maketitle

\newtheorem*{thm}{Theorem}

\theoremstyle{definition}
\newtheorem*{rem}{Remark}
\newtheorem*{ack}{Acknowledgements}

\section{Introduction}
\label{first}
The theorem of Lutz and Martinet~\cite{mart71} asserts that any closed,
oriented $3$-manifold $Y$ admits a contact structure in each homotopy
class of tangent $2$-plane fields. Here $2$-plane
fields are understood to be cooriented; a {\it contact
structure} is understood to be cooriented and positive, that is, a
$2$-plane field $\xi$ defined as the kernel of a global $1$-form $\alpha$
on $Y$ such that $\alpha\wedge d\alpha$ is a positive volume form.

The cited paper by Martinet only covers the existence of {\em some} contact
structure on a given~$Y$; for a proof of the existence of such a structure
in every homotopy class of $2$-plane fields see~\cite{geig}, which gives a
proof of that result along the lines of the original (and never fully
published) argument by Lutz.

The key to that second step is what is nowadays known as a {\it Lutz twist},
a surgery on a knot $T$ in a given contact manifold $(Y,\xi )$ --- with
$T$ {\em transverse} to~$\xi$ --- that is topologically trivial (i.e., does
not change~$Y$), but transforms $\xi$ to a contact structure $\xi'$ in
a different homotopy class of $2$-plane fields.

In a series of papers~\cite{dige01,dige04,dgs04} we described a notion
of contact $r$-surgery on {\it Le\-gendrian} knots in a contact manifold
(that is, knots tangent to the given contact structure), where
$r\in\bfq^*\cup\{\infty\}$ denotes the framing of the surgery relative
to the natural contact framing of the Legendrian knot. This generalises the
contact surgery introduced by Eliashberg~\cite{elia90} and
Weinstein~\cite{wein91}, which in our language is a contact $(-1)$-surgery.
Amongst other things, we
discussed explicit surgery diagrams for various contact manifolds
and gave an alternative proof of the Lutz-Martinet theorem via such
Legendrian surgeries. We did not, however, fully elucidate the relation
between our surgery diagrams and the Lutz twist (although the principal
connection was described in~\cite{dige04}, cf.~\cite{etho01}). The intention
of the present note is to give an explicit Legendrian surgery diagram
for the Lutz twist. In particular, this yields surgery
representations for all contact structures on $S^3$ analogous
to~\cite{dgs04} and provides
concrete realisations for the considerations in Section~6 of~\cite{dige04}.

We shall henceforth assume that the reader is familiar with the basics
of this notion of contact surgery; if not, the best place to start may
well be~\cite{dgs04}, cf.\ also~\cite{ozst}. One fact from~\cite{dige01,dige04}
we should like to recall here is that contact $(+1)$-surgery is the inverse
of contact $(-1)$-surgery.
\section{The Lutz twist}
We briefly recall the definition of the Lutz twist, cf.~\cite{geig}. Let $T$
be a knot transverse to a contact structure $\xi$
on a $3$-manifold $Y$. Then there is a tubular neighbourhood $\nu T$
of $T$ that is contactomorphic to the solid torus $S^1\times D^2_{\delta}$
(with $D^2_{\delta}$ denoting the $2$-disc of radius~$\delta$) for some
suitable $\delta >0$, with contact structure $\zeta =\ker (d\theta +r^2
\, d\varphi )$, where $\theta$ denotes the $S^1$-coordinate, and $r,\varphi$
are polar coordinates on $D^2_{\delta}$. For ease of notation we identify
$(\nu T,\xi)$ with $(S^1\times D^2_{\delta},\zeta)$.

A {\it simple Lutz twist} along $T$ is the operation that replaces the
contact structure $\xi$ on $Y$ by the one that coincides with $\xi$
outside $\nu T$, and on $S^1\times D^2_{\delta}$ is given by
\[ \zeta'=\ker \bigl( h_1(r)\, d\theta +h_2(r)\, d\varphi \bigr),\]
where $h_1,h_2\colon\, [0,\delta ]\rightarrow\bfr$ are smooth functions
satisfying the following conditions:
\begin{itemize}
\item[(i)] $h_1(r)=-1$ and $h_2(r)=-r^2$ for $r$ near $0$,
\item[(ii)] $h_1(r)=1$ and $h_2(r)=r^2$ for $r$ near $\delta$,
\item[(iii)] $(h_1(r),h_2(r))$ is never parallel to $(h_1'(r),h_2'(r))$
  --- in particular, neither of them is ever equal to (0,0) ---,
\item[(iv)] $h_1$ has exactly one zero on the interval $[0,\delta ]$.
\end{itemize}

The boundary conditions (i) and (ii) ensure that $\zeta'$ is defined
around $r=0$ and coincides with $\zeta$ near $r=\delta$; (iii) is the
condition for $\zeta'$ to be a contact structure; condition (iv)
fixes the homotopy class (as $2$-plane field) of the new contact
structure.

The contact structure $\zeta'$ is a so-called {\it overtwisted}
contact structure in the sense of Eliashberg~\cite{elia89},
and as shown in that paper (specifically, Theorem~3.1.1), the
classification of such overtwisted contact structures up to isotopy
fixed near the boundary coincides with the classification of $2$-plane
fields up to homotopy rel boundary. An immediate consequence of that
classification is that the contact structure on $Y$ obtained from $\xi$
by a Lutz twist along $T$ is (up to isotopy) independent of any of the
choices in the construction described above.

\section{The surgery diagram for a Lutz twist}
Let $(Y,\xi )$ be a given contact $3$-manifold and $T$ a knot in $Y$
transverse to $\xi$. In order to describe a Legendrian link
$\Li$ in $Y$ such that $(+1)$-contact surgery on $\Li$ has the same effect as
a Lutz twist along $T$, we may assume by~\cite{dige04} that $(Y,\xi )$
has been obtained from $S^3$ with its standard contact structure
$\xi_{st}$ by contact $(\pm 1)$-surgery on a Legendrian link in
$(S^3,\xi_{st})$, and thus can be represented by the front projection
(to the $yz$-plane) of this Legendrian link, considered as a link
in $\bfr^3$ with its standard contact structure $\xi_{st}=\ker (dz+x\, dy)$,
which is contactomorphic to $(S^3,\xi_{st})$ with a point removed.

\begin{figure}[h]
\centerline{\relabelbox\small
\epsfxsize 5cm \epsfbox{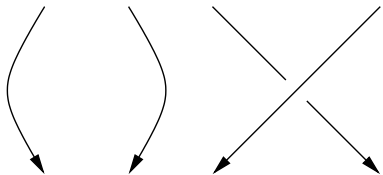}
\endrelabelbox}
\caption{Impossible front projections of positively transverse curve.}
\label{figure:transverse-not}
\end{figure}

For the representation of Legendrian and transverse knots via their
front projection we refer to \cite{efm01,etny,gomp98}.
Beware that these three papers use three different conventions for writing
the standard contact structure on $\bfr^3$. We follow the
one from \cite{gomp98} (which is also that of~\cite{geig}).
The positive transversality condition $\dot{z}+x\dot{y}>0$ for a
curve $t\mapsto (x(t),y(t),z(t))$ implies that in the front projection
of a positively transverse knot there can be no vertical tangencies
going downwards ($\dot{y}=0$, $\dot{z}<0$), and all but the crossing
shown in Figure~\ref{figure:transverse-not} are possible.

From these front projections
it is easy to describe the positive transverse push-off of an oriented
Legendrian knot: smooth the up-cusps and replace the down-cusps
by kinks (there is only one possibility for the sign of the crossing in this
kink). Similarly, one can easily describe an oriented Legendrian knot
whose positive transverse push-off is a given transverse knot,
cf.~\cite{efm01}:
\begin{itemize}
\item[(i)] In the front projection of the given transverse knot (oriented
positively), replace vertical (upwards) tangencies by cusps.
\item[(ii)]  By Figure~\ref{figure:transverse-not}, in those crossings of a
positively transverse knot that cannot be interpreted as the front projection
of a Legendrian knot, at least one of the strands is pointing up $(\dot{z}>0$).
If one adds a zigzag to that strand (if both are going up, either
can be chosen), it is possible to realise the given crossing by the
front projection of a Legendrian curve.
\end{itemize}

Therefore, the following theorem gives a complete surgery description of
Lutz twists.

\begin{figure}[h]
\centerline{\relabelbox\small
\epsfxsize 6cm \epsfbox{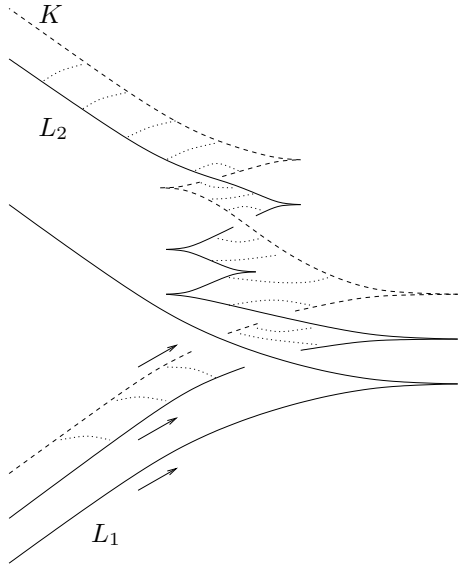}
\extralabel <-4.9cm,0.3cm> {$L_1$}
\extralabel <-5.6cm,5.7cm> {$L_2$}
\extralabel <-5.6cm,7.2cm> {$K$}
\endrelabelbox}
  \caption{Surgery diagram for Lutz twist.}
  \label{figure:lutz}
\end{figure}

\begin{thm}
\label{thm}
Let $L_1$ be an oriented Legendrian knot in $(Y,\xi )$, represented by the
front projection of a Legendrian knot in $(\bfr^3,\xi_{st})$ disjoint from the
link describing $(Y,\xi )$. Let $L_2$ be the Legendrian push-off of $L_1$ with
two additional up-zigzags (see Figure~\ref{figure:lutz}).
Let $\xi'$ be the contact structure on $Y$ obtained from $\xi$ by contact
$(+1)$-surgery on both $L_1$ and $L_2$, and $\xi''$ the contact structure
obtained from $\xi$ by a simple Lutz twist along the positive
transverse push-off $T$ of~$L_1$. Then $\xi'$ and $\xi''$ are isotopic via
an isotopy fixed outside a tubular neighbourhood of~$L_1$.
\end{thm}

The proof of this theorem proceeds as follows: First of all, we verify
that the described surgeries on $L_1$ and $L_2$ taken together do not change
the manifold $Y$. Secondly, we check that the resulting contact structure
is overtwisted by exhibiting an explicit overtwisted disc. Then, again
by Eliashberg's classification of overtwisted contact structures,
and thanks to the fact that the two surgeries only change the contact
structure in a tubular neighbourhood of~$L_1$ (which contains
the overtwisted disc just mentioned), it suffices to show
that the described surgeries and the corresponding Lutz twist have
the same effect on the homotopy class of the contact structure, regarded
as a mere plane field.

\vspace{2mm}

(1) Recall that contact $r$-surgery on a Legendrian knot $L$
means that topo\-logically we perform
surgery with coefficient $r\in\bfq^*\cup\{\infty\}$ relative to
the contact framing of $L$, which is determined by a vector field along
$L$ transverse to the contact structure. In the front projection picture
this corresponds to pushing $L$ in $z$-direction, and it is this what
we mean by {\em the} Legendrian push-off of~$L$. As shown in~\cite{dige01},
contact $(+1)$-surgery along $L$ and contact $(-1)$-surgery along
its Legendrian push-off cancel each other, and in particular do not
change the underlying manifold. Since adding two zigzags to a Legendrian knot
adds two negative twists to its contact framing, we see that topologically
the two contact $(+1)$-surgeries on $L_1$ and $L_2$ are the same
as a contact $(+1)$-surgery along $L_1$ and a $(-1)$-surgery along its
Legendrian push-off, and hence topologically trivial.

We can easily see this directly: Write $t$ for the Thurston-Bennequin invariant
of $L_1$, so that the linking number between $L_1$ and its Legendrian push-off
(or with $L_2$) is given by $\lk (L_1,L_2)=t$. Then the topological
framings (i.e., framings relative to the surface framing) of the surgeries
are $n_1=t+1$ and $n_2=t-1$. After a handle slide
(cf.~\cite[Chapter~5]{gost99}) we may replace the link
$(L_1,L_2)$ by $(L_1,L_2-L_1)$, with linking number $\lk (L_1,L_2-L_1)=
t-n_1=-1$ and framing of $L_2-L_1$ equal to
\[ (L_2-L_1)^2=L_2^2+L_1^2-2\,\lk (L_2,L_1)=n_2+n_1-2t=0.\]
By construction, $L_2-L_1$ is an unknot, and the computation above shows that
it is a $0$-framed meridian of $L_1$, which proves the claim that the
composition of the two surgeries is topologically trivial.

\vspace{2mm}

(2) We next exhibit the overtwisted disc in the manifold obtained by
the contact surgeries along $L_1$ and $L_2$. Let $K$ be the knot
indicated in Figure~\ref{figure:lutz}, i.e., the Legendrian push-off
of $L_1$ with one additional zigzag and with one extra negative
linking with $L_2$. Alternatively, $L_2$ may be regarded as the
Legendrian push-off of $K$ with one additional zigzag. 
The surface framing of $L_2$ determined by the Seifert surface
of the oriented link $(-L_2)\sqcup K$ indicated in Figure~\ref{figure:lutz}
is equal to $t-1$ (the contact framing of~$K$), hence equal to
the topological framing used for the surgery on~$L_2$. So that Seifert
surface glued to the meridional disc used for the surgery on $L_2$
defines a disc with boundary $K$ in the surgered manifold. The
surface framing of $K$ determined by that disc equals the
contact framing $t-1$, which is exactly the condition for
an overtwisted disc.

The above verification that $K$ is the boundary of an overtwisted disc
in the surgered manifold is completely straightforward. Nonetheless,
it may be instructive to see that $K$ is not found by accident.
Start with a meridian $K'$ to
both $L_1$ and $L_2$, that is, an unknot with $\lk (K',L_1)=
\lk (K',L_2)=1$. If surgery along $L_1$ and $L_2$ has any chance of being
a Lutz twist, we expect $K'$ to be isotopic to the
boundary of an overtwisted disc.

There is an obvious pair of pants with boundary the oriented link
$L_1\sqcup (-K')\sqcup (-L_2)$ that gives $K'$ the surface framing $n_{K'}=0$
and $L_1,L_2$ the framing $t+1, t-1$, respectively.
Now perform a handle slide of
$-K'$ over the $2$-handle attached to $L_1$ (corresponding to the surgery)
to form, in the surgered manifold, the knot $L_1-K'$, which is
the knot $K$ from above. We compute the linking
numbers
\begin{eqnarray*}
\lk (L_1-K',L_1) & = & n_1-1\; = \; t,\\
\lk (L_1-K',L_2) & = & \lk (L_1,L_2)-1\; = \; t-1,
\end{eqnarray*}
and the surface framing of $L_1-K'$, now with respect to the
annulus with boundary $(L_1-K')\sqcup (-L_2)$:
\begin{eqnarray*}
(L_1-K')^2 & = & n_1+n_{K'}-2\,\lk (L_1,K')\\
          & = & (t+1)+0-2\cdot 1\\
          & = & t-1,
\end{eqnarray*}
which is exactly what we had found for $K$ before.

\vspace{2mm}

(3) It remains to be shown that the topological effect of the
surgeries described in the theorem has the same effect on
the homotopy class of the contact structure (regarded merely as a plane field)
as a Lutz twist. By the neighbourhood theorem for Legendrian submanifolds
(cf.~\cite{geig}), the particular nature of $L_1$ is irrelevant for
this consideration. It therefore suffices to consider specific
examples for $L_1$, where the effect of the surgeries on the
obstruction classes determining the homotopy type of the plane field
can be computed explicitly.

For the following considerations cf.~\cite{geig}. The tangent bundle of
the solid torus $S^1\times D^2$ being trivial, (cooriented)
tangent $2$-plane fields on $S^1\times D^2$ can be identified with maps
$S^1\times D^2\rightarrow S^2$. Thus, the obstructions to homotopy of
$2$-plane fields on $S^1\times D^2$ rel boundary $T^2$ are in
\[ H^2(S^1\times D^2,T^2;\pi_2(S^2))\cong \bfz\]
and
\[ H^3(S^1\times D^2,T^2;\pi_3(S^2))\cong\bfz .\]
The first obstruction corresponds to the extension of a given $2$-plane field
along $T^2$ over a meridional disc of the solid torus and is detected
by the (relative) first Chern class of the plane field (here
the absence of $2$-torsion is crucial). The second obstruction relates to
the extension of the plane field over the $3$-cell one needs to attach to
$T^2\cup$ (meridional disc) to form the solid torus. This obstruction is
captured by the Hopf invariant.

\vspace{2mm}

(3a) In order to deal with the first obstruction, we consider $Y=S^1\times S^2
\subset S^1\times\bfr^3$ with its standard tight contact structure
$\xi =\ker (x\, d\theta +y\, dz-z\, dy)$, in obvious notation, and take $L_1$
to be an oriented Legendrian knot in the homology class of $S^1\times
\{ pt.\}$. The contact manifold $(S^1\times S^2,\xi )$ can be represented by
contact $(+1)$-surgery on a Legendrian unknot $L_0$ with only two cusps,
see~\cite{dgs04}. For $L_1$ we take another such unknot linked once with
$L_0$, and for $L_2$ its Legendrian push-off with additional zigzags as
in the theorem. Write $\xi'$ for the contact structure on
$Y$  obtained by performing contact $(+1)$-surgery on $L_1$ and~$L_2$.

The contact structure $\xi$ has first Chern class $c_1(\xi )=0$. This
follows from the observation that the vector field
\[ (z-y)\, \partial_x+x\,\partial_y-x\,\partial_z+(y+z)\,\partial_{\theta}\]
defines a trivialisation of~$\xi$. Alternatively, this is a consequence
of the homological computations in Section~3 of~\cite{dgs04}, given the
fact that the rotation number $\rot (L_0)$ (with any
orientation on~$L_0$), which can be computed from the front projection as
$(\# \mbox{\rm (down-cusps)}-\# \mbox{\rm (up-cusps)})/2$, is equal to $0$.

In the sequel we assume that the reader is familiar with those homological
computations. Write $\mu_1,\mu_2$ for the meridional circles to
$L_1,L_2$, respectively, as well as the homology classes they
represent in the homology of the surgered manifold. Then, with
$PD$ denoting the Poincar\'e duality isomorphism from cohomology
to homology,
\begin{eqnarray*}
c_1(\xi') & = & \rot (L_1)PD^{-1}(\mu_1)+\rot (L_2)PD^{-1}(\mu_2)\\
          & = & -2PD^{-1}(\mu_2).
\end{eqnarray*}
(This would be true even if $\rot (L_1)\neq 0$, since $\mu_1+\mu_2$ bounds
a disc in $Y$ also after the surgery.)

Let $L_1'$ be a Legendrian push-off of $L_1$. Then the surgery along
$L_1$ and $L_2$ may be assumed to occur in a tube containing $L_1$
and $L_2$, but not $L_1'$. This implies that $L_1'$
represents the same homology class in $H_1(Y)$ both before and
after the surgery. Since $\lk (L_1',L_2)=t$ and along $L_2$
we perform surgery with topological framing $t-1$, we have that
$L_1'-\mu_2$ is homologically trivial in the surgered manifold. Hence
\[ \mu_2=[L_1']=[L_1]\in H_1(Y),\]
so that
\[ c_1(\xi')=-2PD^{-1}([L_1]),\]
which is the same as for a Lutz twist along the positive transverse 
push-off of $L_1$ (i.e., a transverse knot in the homology class of $L_1$),
see~\cite[Prop.~3.15]{geig}.

Since $[L_1]$ generates $H_1(Y)$ in this example, this fully determines
the effect of the surgery on the $2$-dimensional obstruction class.

\vspace{2mm}

(3b) Finally, in order to see that the effect that the surgery on the link
$\Li = L_1\sqcup L_2$ has on the $3$-dimensional obstruction is the same as
that of a Lutz twist along a positive transverse push-off of $L_1$, it is
sufficient to consider an arbitrary oriented Legendrian knot $L_1$ in
$(S^3,\xi_{st})$. Set $r=\rot (L_1)$, so that $\rot (L_2)=r-2$. As
before we write $t$ for the Thurston-Bennequin invariant of $L_1$,
so that the Thurston-Bennequin invariant of $L_2$ equals $t-2$. Let $X$
be the handlebody obtained from $D^4$ by attaching two $2$-handles
corresponding to the two surgeries. Let $c\in H^2(X)$ be the
cohomology class that evaluates to $\rot (L_i)$ on the surface in $X$ given
by gluing a Seifert surface (with induced orientation) of $L_i$ in $D^4$
with the core disc of the corresponding handle, $i=1,2$. Since we
perform $q=2$ contact $(+1)$-surgeries, Corollary~3.6 of~\cite{dgs04}
tells us that the $3$-dimensional invariant of the contact structure
$\xi'$ obtained by these surgeries is given by
\begin{eqnarray*}
d_3(\xi') & = & \frac{1}{4}\bigl( c^2-3\sigma (X)-2\chi (X)\bigr) +q\\
          & = & \frac{1}{4}c^2-\frac{3}{4}\sigma (X)+\frac{1}{2}.
\end{eqnarray*}
The signature $\sigma (X)$ is the signature of the matrix
$\left(\begin{array}{cc}t+1&t\\t&t-1\end{array}\right)$, hence equal to $0$.
Moreover, by that same formula we have $d_3(\xi_{st})=-1/2$. So the change
in the $d_3$-invariant caused by the surgery is
\[ d_3(\xi')-d_3(\xi_{st})=\frac{1}{4}c^2+1.\]
As shown in Section~3 of~\cite{dgs04}, $c^2$ can be computed as
$ar+b(r-2)$, where $(a,b)$ is the solution of
\[ \left(\begin{array}{cc}t+1&t\\t&t-1\end{array}\right)
\left(\begin{array}{c}a\\b\end{array}\right) =
\left(\begin{array}{c}r\\r-2\end{array}\right) .\]
This yields $a=r-2t$ and $b=2-r+2t$, hence
$c^2=4r-4t-4$ and finally $d_3(\xi')-d_3(\xi_{st})=r-t$. This
is exactly minus the so-called self-linking number $l(T)$ of the positive
transverse push-off $T$ of $L_1$, cf.~\cite{etny}.

As shown in \cite{geig}, the relative $d_3$-invariant $d_3(\xi'',\xi_{st})$,
measuring the obstruction to homotopy over the $3$-skeleton
between $\xi_{st}$ and the contact structure $\xi''$ obtained by a Lutz
twist along $T$, equals $l(T)$. Thus, to conclude the proof one would 
need to verify that the absolute $d_3$-invariant of \cite{gomp98} and the
relative $d_3$-invariant of \cite{geig} are, in the case at hand,
related by
\begin{equation}
d_3(\xi_1,\xi_2)=d_3(\xi_2)-d_3(\xi_1).
\label{eqn:d3}
\end{equation}

This can be done by looking at explicit geometric models, though, as always,
it is difficult to keep track of signs. So here is a more roundabout
algebraic argument. Let $\xi_{\pm 1}$ be the contact structure obtained
from $\xi_{st}$ by a Lutz twist along a transverse knot $T_{\mp 1}$
with self-linking number $l(T_{\mp 1})=\mp 1$ (this sign convention
will be explained below); recall that the self-linking number is
independent of the orientation of the transverse knot.
For any natural number $n$, write $\xi_{\pm n}$ for
the contact structure on $S^3$ given by taking the connected sum of $n$ copies
of $(S^3,\xi_{\pm 1})$. The additivity of the relative $d_3$-invariant
implies $d_3(\xi_{\pm n},\xi_{st})=\mp n$, which means that we get a contact
structure on $S^3$ in each homotopy class of tangent $2$-plane fields.

The absolute $d_3$-invariant --- for $2$-plane fields on~$S^3$ --- takes
all the values in $\bfz +1/2$, with $d_3(\xi_{st})=-1/2$. 
By \cite[Lemma~4.2]{dgs04}, it satisfies the additivity rule
\[ d_3(\eta_1\#\eta_2)=d_3(\eta_1)+d_3(\eta_2)+\frac{1}{2}.\]

These observations imply equation~(\ref{eqn:d3}) up to sign. We conclude
\[ d_3(\xi'',\xi_{st})=l(T)=d_3(\xi_{st})-d_3(\xi')=\pm d_3(\xi',\xi_{st}).\]
By the considerations in (3a),
we know that the extension of the contact structure over a meridional
disc is the same, up to homotopy, for surgery on $\Li$ or Lutz
twist along~$T$. From the fact that there are standard models for the tubular
neighbourhood of a Legendrian or transverse knot, respectively,
we infer that $d_3(\xi',\xi_{st})$ and $d_3(\xi'',\xi_{st})$ can only
differ by a constant term independent of the specific knot (corresponding
to a different extension of the $2$-plane field over the $3$-cell attached to
$T^2\cup$ (meridional disc)). Hence, the equation above can only
hold if that constant is zero and the sign is the positive one.
In turn, this yields equation~(\ref{eqn:d3}) in full generality.

(Our definition of $\xi_{\pm n}$ then entails $d_3(\xi_1)=1/2$ and
$d_3(\xi_{-1})=-3/2$, which accords with our labelling of these
structures in~\cite{dgs04}.)

\vspace{2mm}

This concludes the proof of the theorem.

\begin{rem}
If one defines $L_2$ by adding two down-zigzags instead of up-zigzags,
in (3a) one obtains $c_1(\xi')=2PD^{-1}([L_1])$. This is
the same as for a Lutz twist along the negative transverse push-off $T_-$
of $L_1$, since $T_-$ with the orientation that makes it positively
transverse to $\xi$ represents the class $-[L_1]$.
Similarly, with this $L_2$ we find in (3b) that
$d_3(\xi')-d_3(\xi_{st})$ is equal to minus the self-linking number
$t+r$ of the negative transverse push-off of $L_1$.
Therefore, this choice of $L_2$ amounts to
performing a Lutz twist along~$T_-$.
\end{rem}

\begin{ack}
F.~D.\ is partially supported by grant no.\ 10201003 of the National Natural
Science Foundation of China. H.~G.\ is partially supported by grant no.\
GE 1254/1-1 of the Deutsche Forschungsgemeinschaft within the
framework of the Schwerpunktprogramm 1154 ``Globale
Differential\-geo\-metrie''. A.~S.\ is partially supported by
OTKA T034885.
\end{ack}

\end{document}